\documentclass[reqno,10pt]{amsart}
\usepackage{amsfonts, amsmath, amssymb, amsthm, color}
\allowdisplaybreaks[1]
\catcode`@=11
\newbox\tr@tto
\setbox\tr@tto=\hbox{{\count0=0\dimen0=-,9pt\dimen1=1,1pt\loop\ifnum
\count0<11 \advance \count0 by1 \vrule width.51pt height\dimen1
depth\dimen0\kern-0.17pt\advance\dimen0
by-0.05pt\advance\dimen1by0.1pt\repeat \loop\ifnum\count0<21\advance
\count0 by1 \vrule width.6pt height\dimen1 depth\dimen0\kern-0.2pt
\advance\dimen0 by-0.1pt\advance\dimen1 by 0.05pt\repeat}}
\def\medint{\displaystyle\copy\tr@tto\kern-10.4pt\int}
\catcode`@=12

\newtheorem{thm}{Theorem}[section]
\newtheorem{lemma}{Lemma}[section]

\theoremstyle{definition}
\numberwithin{equation}{section}
\newcommand{\al}{\alpha}

\newcommand{\eps}{\varepsilon}
\newcommand{\la}{\lambda}
\newcommand{\ga}{\gamma}
\newcommand{\Si}{\Sigma}
\newcommand{\wnn}{w_{1,n}}
\newcommand{\wnnn}{w_{2,n}}

\newcommand{\Yun}{Y_{1,n}}
\newcommand{\Ydn}{Y_{2,n}}
\newcommand{\Oo}{\mathcal O}
\newcommand{\calO}{\mathcal O}

\newcommand{\mv}{\kappa}
\newcommand{\twunu}{\widetilde w_{1,n}^1}
\newcommand{\twdnu}{\widetilde w_{2,n}^1}
\newcommand{\twund}{\widetilde w_{1,n}^2}
\newcommand{\twdnd}{\widetilde w_{2,n}^2}
\newcommand{\twuu}{\widetilde w_{1}^1}

\newcommand{\twdd}{\widetilde w_{2}^1}
\newcommand{\Onu}{\mathcal O_n^1}
\newcommand{\Ond}{\mathcal O_n^2}
\newcommand{\tVunu}{\widetilde V_{1,n}^1}
\newcommand{\tVdnu}{\widetilde V_{2,n}^1}
\newcommand{\tVund}{\widetilde V_{1,n}^2}
\newcommand{\tVdnd}{\widetilde V_{2,n}^2}
\newcommand{\Vn}{V_{1,n}}
\newcommand{\Vnn}{V_{2,n}}
\newcommand{\edn}{\varepsilon_{2,n}}
\newcommand{\eun}{\varepsilon_{1,n}}

\newcommand{\Id}{\int_\Sigma e^{-\gamma v}}

\newcommand{\mvn}{\kappa}

\newcommand{\uun}{u_{1,n}}
\newcommand{\udn}{u_{2,n}}

\newcommand{\R}{{\mathbb R}}

\newcommand{\vn}{v_n}

\newcommand{\Iu}{\int_\Sigma e^v}

\newcommand{\dvg}{dV_g}
\newcommand{\lla}{\lambda_1}
\newcommand{\llla}{\lambda_2}
\newcommand{\llan}{\lambda_{1,n}}
\newcommand{\lllan}{\lambda_{2,n}}
\newcommand{\Sp}{\mathcal S_1}
\newcommand{\Sm}{\mathcal S_2}

\newcommand{\mmuu}{\mu_2}
\newcommand{\mmun}{\mu_{1,n}}
\newcommand{\mmuun}{\mu_{2,n}}
\newcommand{\mm}{m_1}
\newcommand{\mmm}{m_2}
\newcommand{\ue}{U_\eps}
\newcommand{\ve}{v_\eps}
\begin{document}
\title[]{Minimal blow-up masses and existence of solutions for an asymmetric sinh-Poisson equation}
\author[T.~Ricciardi]{Tonia Ricciardi}
\address[Tonia Ricciardi] {Dipartimento di Matematica e Applicazioni, 
Universit\`{a} di Napoli Federico II, Via Cintia, Monte S.~Angelo, 80126 Napoli, Italy}
\email{tonricci@unina.it}
\author[G.~Zecca]{Gabriella Zecca}
\address[Gabriella Zecca] {Dipartimento di Matematica e Applicazioni, 
Universit\`{a} di Napoli Federico II, Via Cintia, Monte S.~Angelo, 80126 Napoli, Italy}
\email{g.zecca@unina.it}
\begin{abstract}
For a sinh-Poisson type problem with asymmetric exponents 
of interest in hydrodynamic turbulence, we establish the optimal
lower bounds for the blow-up masses.
We apply this result to construct solutions of mountain pass type on two-dimensional tori.
\end{abstract}
\date{}
\maketitle
%%%%%%%%%%%%%%%%%%%%%%%%%%%%%%%%%%%%%%%%%%%%%%%%%%%%%%%%%%%%%%%%%%%%%%%%%%%%%%%
%%%%%%%%%%%%%%%%%%%%%%%%%%%%%%%%%%%%%%%%%%%%%%%%%%%%%%%%%%%%%%%%%%%%%%%%%%%%%%%
\section{Introduction and main results}
\label{sec:intro}
%%%%%%%%%%%%%%%%%%%%%%%%%%%%%%%%%%%%%%%%%%%%%%%%%%%%%%%%%%%%%%%%%%%%%%%%%%%%%%%
%%%%%%%%%%%%%%%%%%%%%%%%%%%%%%%%%%%%%%%%%%%%%%%%%%%%%%%%%%%%%%%%%%%%%%%%%%%%%%%
Motivated by the equations introduced in the context of hydrodynamic turbulence by Onsager in \cite{On} 
and by Sawada and Suzuki in \cite{SawadaSuzuki}, 
we study the following sinh-Poisson type problem:
\begin{equation}
\label{eq:mainpb}
\left\{
\begin{aligned}
-\Delta_g v=&\lla\frac{e^v}{\int_\Si e^v\,\dvg}-\llla\frac{e^{-\ga v}}{\int_\Si e^{-\ga v}\,\dvg}-\kappa&&\mbox{on\ }\Si\\
\int_\Si v\,\dvg=&0.
\end{aligned}
\right.
\end{equation}
Here, the unknown function $v$ corresponds to the stream function of the turbulent Euler flow, 
$(\Si,g)$ is a compact, orientable, Riemannian 2-manifold  without boundary,  the constant 
$\ga\in(0,1]$ describes the intensity of the negatively oriented vortices,
and $\la_1,\llla>0$ are constants related to the inverse temperature, $\mv\in\mathbb R$ is given by
\[
\kappa=\frac{\lla-\llla}{|\Si|},
\]
where $|\Sigma|=\int_\Sigma\dvg$ is  the volume of $\Sigma$. The constant $\mv$ ensures that the right hand side in \eqref{eq:mainpb}
has zero mean value. The trivial solution $v\equiv0$ always exists for \eqref{eq:mainpb}.
\par
Problem \eqref{eq:mainpb} admits a variational structure. Indeed,  \eqref{eq:mainpb} is the Euler-Lagrange equation 
of the functional:
\[
J(v)=\frac12\int_\Sigma |\nabla v|^2\,\dvg -\lambda_1 \log \int_\Sigma e^v\dvg 
-\frac{\llla}{\gamma} \log\int_\Sigma e^{-\gamma v}\,\dvg
\]
defined on the space 
\[
\mathcal E=\left\{v\in H^1(\Sigma):\ \int_\Sigma v\,\dvg=0\right\}. 
\]
With this notation, $v$ is a classical solution for \eqref{eq:mainpb} if and only if
$v$ is a critical point for $J$ in $\mathcal E$.
\par
When $\ga=1$, problem~\eqref{eq:mainpb} reduces to the mean field equation for equilibrium turbulence
derived in \cite{JoyceMontgomery, PointinLundgren}, under the assumption that the point vortices
have equal intensities and arbitrary orientation.
This case has received a considerable attention in recent years. In particular, existence of saddle point
solutions was obtained in \cite{Jevnikar, Ric2007, Zhou}; such results exploit the blow-up analysis
derived in \cite{JostWangYeZhou, OhtsukaSuzuki2006} in order to establish the compactness of solution sequences.
A detailed blow-up analysis is also contained in \cite{JevnikarWeiYang}.
By the Lyapunov-Schmidt approach introduced in \cite{EspositoGrossiPistoia},
sign-changing concentrating solutions for a related local sinh-Poisson type equation were constructed in \cite{BartolucciPistoia} and 
sign-changing bubble-tower solutions were constructed in \cite{GrossiPistoia}.
\par
The case $\ga\neq1$ corresponds to the assumption where all positively oriented point vortices have unit intensity
and all negatively oriented point vortices have intensity equal to $\ga$. 
Such a case already appears in an unpublished manuscript of Onsager \cite{EyinkSreenivasan},
although the rigorous derivation of \eqref{eq:mainpb} is due to \cite{SawadaSuzuki}. 
Actually, the mean field equation derived in \cite{SawadaSuzuki} concerns the very general case where the point vortex intensity
is described by a Borel probability measure $\mathcal P\in\mathcal M([-1,1])$;
it is given by
\begin{equation}
\label{eq:SS}
\left\{
\begin{aligned}
-\Delta_g v=&\la\int_{[-1,1]}\al\left(\frac{e^{\al v}}{\int_\Sigma e^{\al v}\,\dvg}
-\frac{1}{|\Sigma|}\right)\,\mathcal P(d\al)&&\mbox{on\ }\Si\\
\int_\Si v\,\dvg=&0,
\end{aligned}
\right.
\end{equation}
where $\la>0$ is a constant related to the inverse temperature.
The variational functional for \eqref{eq:SS} takes the form
\begin{equation}
\label{def:JP}
J_{\mathcal P}(v)=\frac12\int_\Sigma|\nabla v|^2\,\dvg
-\la\int_{[-1,1]}\log\left(\int_\Sigma e^{\al v}\,\dvg\right)\,\mathcal P(d\al), \qquad v\in\mathcal E.
\end{equation}
In this context, problem~\eqref{eq:mainpb} corresponds to the case
\begin{equation}
\label{def:specialP}
\mathcal P=\tau\delta_1+(1-\tau)\delta_{-\gamma}, \qquad \la_1=\la\tau,\ \la_2=\la\ga(1-\tau).
\end{equation}
Problem~\eqref{eq:SS} was analyzed in \cite{ORS, RZ}, where concentration properties 
for non-compact solution sequences is studied in the spirit of \cite{BrezisMerle}.
The optimal Moser-Trudinger constant for $J_\mathcal P$
is computed in \cite{RS}.
More precisely, it is shown in \cite{RS} that $J_{\mathcal P}$ is bounded from below
if $\lambda<\bar\lambda$, where
\begin{equation}
\label{def:barlambda}
\bar\lambda=8\pi\inf\left\{\frac{\mathcal P(K_\pm)}{\left(\int_{K_\pm}\al\,\mathcal P(d\al)\right)^2}:\ 
K_\pm\subset I_\pm\cap\mathrm{supp}\mathcal P\right\},
\end{equation}
where $I_+:=[0,1]$ and $I_-=[-1,0)$, $K_\pm$ are Borel sets. On the other hand, 
$J_{\mathcal P}$ is unbounded from below if $\la>\bar\lambda$.
In the case where $\mathcal P$ is discrete, an equivalent form of \eqref{def:barlambda} was derived
in \cite{ShafrirWolansky} in the context of Liouville systems.
However, particularly with respect to the blow-up behavior of solutions, 
it appears complicated to further analyze \eqref{eq:SS} in its full generality.
Therefore, recent effort has been devoted to study cases where $\mathcal P$
is the sum of two Dirac masses, such as \eqref{eq:mainpb}.
In particular, for \textit{small} values of $\gamma$, existence results for \eqref{eq:mainpb} in planar domains with holes  
were obtained in \cite{RTZZ}; 
sign-changing blow-up solutions for a local Dirichlet problem related to \eqref{eq:mainpb}
are constructed in \cite{PistoiaRicciardi}.
The case $\ga=1/2$ is equivalent to the Tzitz\'eica equation describing surfaces of constant affine curvature,
and was recently analyzed in \cite{JevnikarYang2016}.
\par
Here, we are concerned with the existence and blow-up properties of \eqref{eq:mainpb}
for \textit{general} values of $\ga\in(0,1]$.
It is known \cite{BrezisMerle, ORS, RZ} that unbounded solution sequences for problem~\eqref{eq:mainpb} necessarily concentrate
on a finite set $\mathcal S\subset\Sigma$.
Our first aim in this note is to derive the optimal lower bounds for the blow-up masses,
see Theorem~\ref{thm:blowup} below.
\par
Then, we apply the blow-up analysis results to the construction of a non-zero solution for \eqref{eq:mainpb}
via a mountain pass argument, provided the compact surface $\Sigma$
satisfies a suitable condition concerning the first eigenvalue, see Theorem~\ref{thm:mp} below.
%\par
%Problem~\eqref{eq:mainpb} always admits the trivial solution $v\equiv0$. 
%Our aim in this note is to study the existence of nontrivial solutions to \eqref{eq:mainpb} for suitable values of $\lambda_1$ and $\lambda_\ga$. 
%It is well-known that Problem~\eqref{eq:mainpb} presents a lack of compactness.
%More precisely, along a solution sequence $(\llan,\lllan,\vn)$, the following measures may concentrate up to subsequences,
%namely:
%\begin{equation}
%\llan\frac{e^{\vn}}{\int_\Si e^{\vn}\,dx}\stackrel{\ast}{\rightharpoonup}\sum_{x_0\in\mathcal S_+}m_+(x_0)\delta_{x_0}+r_+(x),
%\qquad
%\lllan\frac{e^{-\ga\vn}}{\int_\Si e^{-\ga\vn}\,dx}\stackrel{\ast}{\rightharpoonup}\sum_{y_0\in\mathcal S_-}m_-(y_0)\delta_{y_0}+r_-(x),
%\end{equation}
%Our first aim is to refine the known existence results for $m_\pm$ and $r_\pm$, in particular by estimating the sharp minimal
%blow-up masses.
%Then, we apply such results to prove the existence of mountain pass type solutions.
\par
In order to state our results precisely, we introduce some notation.
Let $(\llan,\lllan,\vn)$ be a solution sequence to \eqref{eq:mainpb}.
If $\vn$ is unbounded in $L^\infty(\Sigma)$,
then, up to subsequences, we have
\begin{equation}
\label{eq:concentration}
\begin{aligned}
&\llan\frac{e^{\vn}}{\int_\Si e^{\vn}\,\dvg}\,\dvg\stackrel{\ast}{\rightharpoonup}\sum_{p\in\Sp}m_1(p)\delta_{p}+r_1\,\dvg,\\
&\lllan\frac{e^{-\ga\vn}}{\int_\Si e^{-\ga\vn}\,\dvg}\,\dvg\stackrel{\ast}{\rightharpoonup}\sum_{p\in\Sm}m_2(p)\delta_{p}+r_2\,\dvg,
\end{aligned}
\end{equation}
weakly in the sense of measures, where the blow-up sets $\Sp,\Sm\subset\Si$ are finite, 
the \lq\lq blow-up masses" $m_i(p)$ satisfy the lower bound $m_i(p)\ge4\pi$ for all
$p\in\mathcal S_1\cup\mathcal S_2$, 
$r_i\in L^1(\Si)$, $i=1,2$.
See Lemma~\ref{lem:BM} below for a more detailed statement.
Our first aim is to improve the lower bound for the blow-up masses $m_i$, $i=1,2$.
\begin{thm}
\label{thm:blowup}
Let $(\llan,\lllan,\vn)$ be a concentrating solution sequence for \eqref{eq:mainpb}
and suppose that \eqref{eq:concentration} holds. 
Then, 
the blow-up masses satisfy the following lower bounds:
\begin{equation}
\label{eq:improvedmasses}
\begin{aligned}
&m_1(p)\ge8\pi\ \forall p\in\Sp;
&&m_2(p)\ge\frac{8\pi}{\ga}
\ \forall p\in\Sm.
\end{aligned}
\end{equation}
%\item[(ii)]
%Suppose that $\mathcal S_i\setminus\mathcal S_j\neq\emptyset$
%for a fixed $i\in\{1,2\}$, and $j\neq i$.
%Then, $r_i\equiv0$.
\end{thm}
We apply Theorem~\ref{thm:blowup} to derive the existence of non-zero
mountain-pass solutions to problem~\eqref{eq:mainpb}.
In order to state the existence result,
we denote by $\mu_1(\Si)$ the first positive eigenvalue of $-\Delta_g$ on $\Sigma$, namely
\begin{equation}
\label{def:mu1}
\mu _{1}(\Sigma):=\inf_{\phi\in \mathcal E\setminus \{0\}}\frac{\int_{\Sigma
}\left\vert \nabla  \phi \right\vert ^{2}\dvg}{\int_{\Sigma }\phi ^{2} \,\dvg}.
\end{equation}
Our existence result holds for surfaces $\Sigma$ satisfying the condition
\begin{equation}
\label{ipmu}
8\pi<\mu_1(\Sigma)|\Sigma|< 16\pi(1+\ga).
\end{equation}
Condition~\eqref{ipmu} is satisfied, e.g., when $\Sigma$ is  the flat torus $\R^2/\mathbb Z^2$, 
since  in this case $\mu_1(\Sigma)|\Sigma|=4\pi^2$.
Let
\begin{equation}
\label{def:Lambda}
\Lambda:=\left\{
(\lla,\llla)\subset\mathbb R^2: 
\begin{aligned}
&(i)\quad &&\la_1,\llla\ge0, \max\{ \la_1,\gamma\llla\} >8\pi \\
&(ii)\quad &&\la_1\not\in 8\pi\mathbb N, \ \llla\not\in\frac{8\pi}{\ga}\mathbb N\\
&(ii)\quad\ &&\la_1+\ga\llla<\mu_1(\Sigma)|\Sigma|
\end{aligned}
\right\}.
\end{equation}
Our existence result is the following.
\begin{thm}
\label{thm:mp}
Suppose $(\Si,g)$ satisfies \eqref{ipmu}.
Then, there exists a nontrivial solution to \eqref{eq:mainpb} for all $(\la_1,\llla)\in\Lambda$.
\end{thm}
We note that condition \eqref{ipmu} implies $\Lambda\neq\emptyset$. 
More precisely,
$\Lambda$ is a union of two (possibly overlapping) triangles $T_1$, $T_2$, where  
$T_1$ has vertices $\{(8\pi,0),\,(\mu_1(\Sigma)|\Sigma| ,0),\,(8\pi, (\mu_1(\Sigma)|\Sigma|-8\pi)/\ga)\}$, 
and $T_2$ 
has vertices $\{(0,8\pi/\ga),\,(0, \mu_1(\Sigma)|\Sigma|/\ga),\,( \mu_1(\Sigma)|\Sigma|-8\pi, 8\pi/\ga)\}$.
\par
%%%%%%%%%%%%%%%%%%%
%%%%%%%%%%%%%%%%%%%
The remaining part of this note is devoted to the proof of Theorem~\ref{thm:blowup}
and Theorem~\ref{thm:mp}.
The proof of Theorem~\ref{thm:blowup} is in the spirit of \cite{OhtsukaSuzukiBanachCenter, RZmathJap+uno},
where the case $\ga=1$ is considered. The case $\ga\neq1$ introduces an asymmetry in the problem, 
which requires some careful modifications in the proof.
The proof of Theorem~\ref{thm:mp} relies on the variational setting introduced in
\cite{StTa, LN, Ric2007}, based on Struwe's monotonicity trick \cite{Struwe}.
In view of the improved lower bounds for blow-up masses, as stated in Theorem~\ref{thm:blowup},
Theorem~\ref{thm:mp} with $\ga=1$ is more general than the corresponding result in
\cite{Ric2007}. 
\subsection*{Notation}
We denote by $C>0$ a general constant whose actual value may vary from line to line.
When the integration variable is clear, 
we may omit it. We take subsequences without further notice.
We set $\|\cdot\|=\|\cdot\|_{\mathcal E}$.
%%%%%%%%%%%%%%%%%%%%%%%%%%%%%%%%%%%%%%%%%%%%%%%%%%%%%%%%%%%%%%%%%%%%%%%%%%%%%%%%%%%%%%%%%%%%%%%%%%%%%%%%%%%%%%%%%%%%%%%%%%%%%%%%%%%%%
%%%%%%%%%%%%%%%%%%%%%%%%%%%%%%%%%%%%%%%%%%%%%%%%%%%%%%%%%%%%%%%%%%%%%%%%%%%%%%%%%%%%%%%%%%%%%%%%%%%%%%%%%%%%%%%%%%%%%%%%%%%%%%%%%%%%%
\section{Proof of Theorem~\ref{thm:blowup}}
\label{sec:blowup}
%%%%%%%%%%%%%%%%%%%%%%%%%%%%%%%%%%%%%%%%%%%%%%%%%%%%%%%%%%%%%%%%%%%%%%%%%%%%%%%%%%%%%%%%%%%%%%%%%%%%%%%%%%%%%%%%%%%%%%%%%%%%%%%%%%%%%
%%%%%%%%%%%%%%%%%%%%%%%%%%%%%%%%%%%%%%%%%%%%%%%%%%%%%%%%%%%%%%%%%%%%%%%%%%%%%%%%%%%%%%%%%%%%%%%%%%%%%%%%%%%%%%%%%%%%%%%%%%%%%%%%%%%%%
This section is devoted to the proof of Theorem~\ref{thm:blowup}.
We denote by $(\llan,\lllan,\vn)$ a solution sequence to \eqref{eq:mainpb}
with $\llan\to\lambda_{1,0}\in\mathbb R$, $\lllan\to\lambda_{2,0}\in\mathbb R$,
as $n\to\infty$.
We define the measures $\mu_{1,n},\mmuun\in\mathcal M(\Sigma )$ by
\begin{align}
\label{mui}
&\mmun=\llan\frac{e^{v_n}}{ \int_\Sigma e^{v_n}}\,\dvg,
&&\mmuun=\lllan\frac{e^{-\gamma \vn}}{\int_\Sigma e^{-\gamma \vn}}\,\dvg.
\end{align}
We may assume that
$\mu_{i,n}\stackrel{*}{\rightharpoonup}\mu_i\in\mathcal M(\Sigma)$ weakly in the sense of measures,
$i=1,2$.
We define the blow-up sets:
\begin{align*}
\Sp=&\{ p\in\Sigma :\exists\,p_n\rightarrow p\ \mathrm{s.t.\ }  v_n(p_n) \rightarrow +\infty)\}\\
\Sm=&\{ p\in\Sigma :\exists\,p_n\rightarrow p\ \mathrm{s.t.\ }  v_n(p_n) \rightarrow -\infty)\}
\end{align*}
and we denote $\mathcal S=\Sp\cup\Sm$.
\par
For the reader's convenience and in order to fix notation, we collect in the following lemma the 
necessary known blow-up results for \eqref{eq:mainpb}.
\begin{lemma}[\cite{ORS,RZ}]
\label{lem:BM}
For the solution sequence $ (\llan,\lllan,\vn)$ the following alternative holds.
\begin{enumerate}
\item[1)] Compactness: $\limsup_{n\rightarrow \infty}\|v_n\|_\infty<+\infty$. 
We have $\mathcal S=\emptyset$
and there exist a  solution $v \in\mathcal E$ to \eqref{eq:mainpb} with  $\lambda_{1}=\lambda_{1,0}$ and $\llla=\lambda_{2,0}$
and a subsequence $\{v_{n_k}\}$ such that $v_{n_k}\rightarrow v$
in $\mathcal E$.
\item[2)] Concentration: $\limsup_{n\rightarrow \infty}\|v_n\|_\infty=+\infty$. We have $\mathcal S\neq\emptyset $ and
    \begin{align*}
    &\mu_1=\sum_{p\in\Sp}\mm(p)\delta_p+r_1\,\dvg\\
&\mmuu=\sum_{p\in\Sm}\mmm(p)\delta_p+r_2\,\dvg
\end{align*}
where $\delta_p$ denotes the Dirac delta centered at $p\in\mathcal S$. 
The constants $\mm(p),\,\mmm(p)$
satisfy the lower bound
\begin{equation}
\label{8pi}
m_i(p)\geqslant 4 \pi,\qquad i=1,2
    \end{equation}
and $r_i\in L^1(\Sigma)\cap L^\infty_{loc}(\Sigma\setminus \mathcal S)$.
Moreover,
for every $p\in\Sp\cap\Sm$, we have
\begin{equation}
\label{identquad}
8\pi\left[\mm(p)+\frac{\mmm(p)}{\gamma}\right]= [\mm(p)-\mmm(p)]^2.
\end{equation}
\end{enumerate}
\end{lemma}
We note that \eqref{identquad} implies that for all
$p\in \Sp\setminus\Sm$ we have $\mmm(p)=0$ and therefore $\mm(p)=8\pi$. Similarly,
for every  $p\in\Sm\setminus\Sp$ we have $\mm(p)=0$
and therefore $\mmm(p)=\frac{8\pi}\ga$.
In particular, Theorem~\ref{thm:blowup} is satisfied for 
$p\in(\Sp\setminus\Sm)\cup(\Sm\setminus\Sp)$.
Therefore, henceforth we assume $p\in\Sp\cap\Sm$.
Let $\calO\subset\mathbb R^2$ be a smooth bounded open set.
We shall need the following version of the Brezis-Merle type alternatives for Liouville systems.
\begin{lemma}[\cite{LN}]
\label{lem:BMsystem}
Suppose $(w_{1,n},w_{2,n})$ is a solution sequence to
the Liouville system
\begin{equation}
\label{Liouvillen2}
\left\{
\begin{split}
&-\Delta w_{1,n}=aV_{1,n}e^{w_{1,n}}-bV_{2,n}e^{ w_{2,n}}&&\mathrm{in\ }\calO\\
&-\Delta w_{2,n}=-cV_{1,n}e^{w_{1,n}}+dV_{2,n}e^{w_{2,n}}&&\mathrm{in\ }\calO,
\end{split}\right.
\end{equation}
where $V_{i,n}\in L^\infty(\calO)$, $i=1,2$, are given functions,
$a,b,c,d>0$ are fixed constants
and
\begin{align}
\label{Vassumpt}
0\leqslant V_{i,n}\le C,\quad \int_\calO e^{w_{i,n}}\le C,\qquad i=1,2.
\end{align}
Then, up to subsequences, exactly one of the following alternatives holds true.
\begin{item}
\item{1.} Both $w_{1,n}$ and $w_{2,n}$ are locally uniformly bounded in $\calO$.
\item{2.} There is  $i\in\{1,2\}$ such that $w_{i,n}$ is uniformly bounded in $\calO$
and $w_{j,n}\rightarrow-\infty $ locally uniformly  in $\calO$ for $j\neq i.$
\item{3.} Both $w_{1,n}\rightarrow -\infty$ and $w_{2,n}\rightarrow -\infty$ locally uniformly in $\calO$.
\item{4.} For the blow-up sets $\mathcal S_1^0$, $\mathcal S_2^0$ defined for this subsequence,
we have $\mathcal S_1^0\cup\mathcal S_2^0\neq\emptyset$ and $\sharp(\mathcal S_1^0\cup\mathcal S_2^0)< +\infty.$
Furthermore, for each $i\in\{1,2\}$, either $w_{i,n}$ is locally uniformly bounded in $\calO\setminus (S_1^0\cup\mathcal S_2^0)$ or
$w_{i,n}\rightarrow -\infty $ locally uniformly  in $\calO\setminus (S_1^0\cup\mathcal S_2^0).$
Here, if $\mathcal S_i^0\setminus (S_1^0\cap\mathcal S_2^0) \neq \emptyset$ then $w_{i,n}\rightarrow -\infty $
locally uniformly  in $\calO\setminus (S_1^0\cup\mathcal S_2^0),$ and each $x_0 \in \mathcal S_i^0$ satisfies 
$m_i(x_0 )\ge4\pi$ such that
\[
V_{i,n}(x)e^{w_{i,n}}\rightharpoonup \sum_{x_0\in\mathcal S_i^0}m_i(x_0) \delta_{x_0} \qquad *\mbox{-weakly in }\mathcal M(\calO),
\]
$i=1,2$.
\end{item}
\end{lemma}
\begin{proof}
For the case $a=c=2$, $b=d=1$, corresponding to the case
of the Toda system, Lemma~\ref{lem:BM} was established in \cite{LN}, Theorem~4.2.
However, by carefully inspecting the proof in \cite{LN}, it is clear that
Theorem~4.2 in \cite{LN} holds true for general $a,b,c,d>0$. 
\end{proof}
We also use the following known result.
\begin{lemma}
\label{lemmaprel}
Let $p_0\in\mathcal S_1\cap\mathcal S_2 $. There exists a sequence $x_{1,n}\rightarrow p_0$ and a sequence $x_{2,n} \rightarrow p_{0}$ such that:
\begin{item}[ ]{i)} $v_n(x_{1,n})\rightarrow +\infty,\qquad\qquad  v_n(x_{1,n})-\log\int_\Sigma e^{v_n}\rightarrow +\infty$,
\item{ii)} $-v_n(x_{2,n})\rightarrow +\infty,\qquad\quad - v_n(x_{2,n})-\frac 1\gamma \log\int_\Sigma e^
{-\gamma v_n}\rightarrow +\infty.$
\end{item}
\end{lemma}
We can now complete the proof of our first result.
\begin{proof}[Proof of Theorem~\ref{thm:blowup}]
We denote by $(\Psi, \mathcal U)$ an isothermal chart satisfying $\bar{\mathcal U}\cap\mathcal S=\{p_0\}$,
$\Psi(\mathcal U)=\calO\subset\mathbb R^2$ and
\begin{equation}
\label{isot}
\Psi(p_0)=0,\qquad g(X)=e^{\xi(X)}(dX_1^2+dX^2_2), \qquad\xi(0)=0,
\end{equation}
where $X=(X_1,X_2)$ denotes a coordinate system on $\calO$.
In particular, identifying $\varphi(X)=v(\Psi^{-1}(X))$ for any function $\varphi$ defined on $\Sigma$,
we have that a solution $v$ to equation \eqref{eq:mainpb} satisfies
\[
-\Delta v=\lambda_1e^\xi\frac{e^v}{\Iu}-\lambda_2e^\xi\frac{e^{-\gamma v}}{\Id}-\mv e^\xi
\qquad\mathrm{in\ }\calO.
\]
We define $h_\xi$ by
\begin{equation}
\label{hxi}
\left\{
\begin{split}
-\Delta h_\xi=&e^{\xi}\qquad\mathrm{in\ }\calO,\\
h_\xi=&0\qquad\mathrm{on\ }\partial\calO.
\end{split}
\right.
\end{equation}
For the solution sequence $ (\la_{1,n},\la_{2,n},v_n)$ to \eqref{eq:mainpb} such that $\la_{1,n}\to\la_{1,0}$ and $\la_{2,n}\to\la_{2,0}$, 
we define $w_{i,n}:\Oo\to\mathbb R$, $i=1,2$, by setting
\begin{align*}
&w_{1,n}:=v_n-\log \int_\Sigma e^{v_n} +\mvn h_\xi,\\
&w_{2,n}:=-\ga v_n-\log\int_\Sigma e^{-\gamma v_n}-\ga\mvn h_\xi,
\end{align*}
where, as before, we identify $v_n(X)=v_n(\Psi^{-1}(X))$.
Then, the above definitions imply that
\begin{align*}
&\frac{e^{\vn}}{\int_\Sigma e^{\vn}}=e^{\wnn-\mvn h_\xi},
&&\frac{e^{-\ga\vn}}{\int_\Sigma e^{-\ga\vn}}=e^{\wnnn+\ga\mvn h_\xi}.
\end{align*}
Setting
\begin{align}
\label{def:V}
&\Vn:=\lambda_{1,n} e^{\xi-\mv h_\xi},
&&\Vnn:=\la_{2,n}e^{\xi+\ga\mv h_\xi},
\end{align}
we may write
\begin{equation}
\label{eq:masses}
\begin{aligned}
\Vn e^{\wnn}=&\lambda_{1,n} e^{\xi-\mv h_\xi} \frac{e^{\vn}}{\int_\Sigma e^{\vn}}e^{\mv h_\xi}
=\lambda_{1,n} e^\xi \frac{e^v}{\Iu}\\
\Vnn e^{\wnnn}=&\la_{2,n}e^{\xi+\ga\mv h_\xi}\frac{e^{-\ga\vn}}{\int_\Sigma e^{-\ga\vn}}e^{-\ga\mvn h_\xi}
=\la_{2,n}e^\xi\frac{e^{-\ga\vn}}{\int_\Sigma e^{-\ga\vn}}
\end{aligned}
\end{equation}
By definition of $w_{1,n}$, $w_{2,n}$ and $h_\xi$, we thus have
\begin{equation*}
-\Delta w_{1,2}=-\Delta\vn+\mv e^\xi
=\lambda_{1,n} e^\xi \frac{e^v}{\Iu}-\la_{2,n}e^\xi\frac{e^{-\ga\vn}}{\int_\Sigma e^{-\ga\vn}}
\end{equation*}
and
\begin{equation*}
\Delta w_{2,n}=\ga\Delta\vn-\ga\mv e^\xi
=-\ga\lambda_{1,n} e^\xi \frac{e^v}{\Iu}+\ga\la_{2,n}e^\xi\frac{e^{-\ga\vn}}{\int_\Sigma e^{-\ga\vn}}.
\end{equation*}
In particular, $(\wnn, \wnnn)$ is a solution to the Liouville system
\begin{equation}
\label{Liouvillen}
\left\{
\begin{split}
&-\Delta w_{1,n}=V_{1,n}e^{w_{1,n}}-V_{2,n}e^{\gamma w_{2,n}}&&\mathrm{in\ }\calO\\
&-\Delta w_{2,n}=-\ga V_{1,n}e^{w_{1,n}}+\ga V_{2,n}e^{\gamma w_{2,n}}&&\mathrm{in\ }\calO,
\end{split}\right.
\end{equation}
which is of the form~\eqref{Liouvillen2} with $a=b=1$, $c=d=\ga$.
We claim that solutions to \eqref{Liouvillen} satisfy the estimates \eqref{Vassumpt} 
in view of \eqref{eq:masses}.
Indeed, we have
\[
0\le V_{i,n}\le\sup_n\la_{i,n}e^{\|\xi\|_\infty+\mv\|h_\xi\|_\infty}\le C,\qquad i=1,2.
\]
Moreover,
\begin{align*}
&\int_{\calO}e^{w_{1,n}}=\int_{\calO}\frac{e^{\vn}}{\int_\Sigma e^{\vn}}e^{\mv h_\xi}
\le e^{\mv\|h_\xi\|_\infty}\\
&\int_{\calO} e^{w_{2,n}}=\int_{\calO}\frac{e^{-\ga\vn}}{\int_\Sigma e^{-\ga\vn}}e^{-\ga\mv h_\xi}
\le e^{\ga\mv\|h_\xi\|_\infty}.
\end{align*}
\par
For later use, we also note that
\begin{equation}
\label{wid}
w_{1,n}+\frac{1}{\ga}w_{2,n}=-\log\Iu- \frac 1\ga\log\Id \leqslant C.
\end{equation}
We define
\begin{equation*}
\begin{split}
V_1=\lambda_{1,0}\,e^{\xi-\mv h_\xi},
\qquad V_2=\lambda_{2,0}\,e^{\xi+\ga\mv h_\xi},
\end{split}
\end{equation*}
so that $V_{i,n}\to V_i$, $i=1,2$,
uniformly on $\overline\calO$. We also define
\[
\mathcal S_i^0= \left\{ X\in \mathcal O : \exists X_n\rightarrow X\ \mathrm{s.t.\ } w_{i,n}(X_n)\rightarrow +\infty\right\},
\quad i=1,2.
\]
In view of Lemma~\ref{lemmaprel} there exist $x_{1,n}$ and $x_{2,n}$ such that
$x_{1,n}\rightarrow p_0$, $v_n(x_{1,n})\rightarrow +\infty$, $x_{2,n}\rightarrow p_0$ and $-v_n(x_{2,n})\rightarrow +\infty$,
%\begin{equation*}
%\begin{array}{lll}
%x_{1,n}\rightarrow p_0& \mbox{ and }& v_n(x_{1,n})\rightarrow +\infty\\
%x_{2,n}\rightarrow p_0& \mbox{ and } & -v_n(x_{2,n})\rightarrow +\infty\\
%\end{array}
%\end{equation*}
and furthermore
\begin{equation*}
\begin{array}{lll}
X_{1,n}= \Psi (x_{1,n})\rightarrow 0& \mbox{ and }& w_{1,n}(X_{1,n})\rightarrow +\infty\\
X_{2,n}= \Psi (x_{2,n})\rightarrow 0& \mbox{ and }& w_{2,n}(X_{2,n})\rightarrow +\infty.
\end{array}
\end{equation*}
In particular, $0\in \mathcal S_1^0\cap \mathcal S_2^0$ and
\begin{equation}
\label{Sio}
\mathcal S_1^0=\Psi(\mathcal U\cap\mathcal S_1)=\left\{0\right\}=\Psi(\mathcal U\cap\mathcal S_2)
=\mathcal S_2^0.
\end{equation}
On the other hand, in view of \eqref{eq:masses} and
recalling that $\xi(0)=0$, we note that
\begin{equation}
\begin{split}
&V_{1,n}e^{w_{1,n}}\,dX\stackrel{*}\rightharpoonup m_1(p_0)\delta_{X=0}+s_1(X)\,dX,\\
&V_{2,n}e^{w_{2,n}}\,dX\stackrel{*}\rightharpoonup m_2(p_0)\delta_{X=0}+s_2(X)\,dX,\\
\end{split}
\end{equation}
$\ast$-weakly in $\mathcal M (\bar{\mathcal O)}$, where $s_i(X)=r_i(\Psi^{-1}(X))e^{\xi(X)}$,
$s_i\in L^1(\mathcal O)\cap L^\infty_{loc}(\bar{\mathcal O} \setminus \{0\})$
and $m_i(p_0)\ge4\pi$, $i=1,2$.
In view of \eqref{Sio} there exist $\Yun,\Ydn\in\Oo$,
$\Yun,\Ydn\to0$ such that
\begin{align*}
w_{i,n}(Y_{i,n})=&\sup_{\mathcal O}w_{i,n}\rightarrow +\infty, \quad i=1,2.
\end{align*}
In order to conclude the proof, we rescale the Liouville system~\eqref{Liouvillen} twice. 
%%%%%%%%%%%%%%%%%%%%%%%%%%%%%%%%%%%%%%%%%%%%%%%%%%%%%%%%%%%%%%%%%%%%%%%%%%%%%%%%%%%%% 
\par
\textit{Proof of $m_1(p_0)\ge8\pi$.}
We first rescale \eqref{Liouvillen} around $\Yun$ with respect the rescaling parameter
\begin{align*}
\eun=&e^{-w_{1,n}(\Yun)/2},
\end{align*}
that is $w_{1,n}(\Yun)=\sup_{\calO}w_{1,n}=-2\log\eun$.
\par
Namely, we define the expanding domain
\[
\Onu=\{X\in\mathbb R^2:\ \Yun+\eun X\in\calO\}
\]
and we define $\twunu,\twdnu:\Onu\to\mathbb R$ by setting
\begin{align*}
\twunu(X)=&w_{1,n}(\Yun+\eun X)-w_{1,n}(\Yun)\\
\twdnu(X)=&w_{2,n}(\Yun+\eun X)-w_{1,n}(\Yun).
\end{align*}
Then, $\twunu,\twdnu$ is a solution for the Liouville system
\begin{equation}
\label{Liouville1}
\left\{
\begin{split}
-\Delta\twunu=&\tVunu e^{\twunu}-\tVdnu e^{\twdnu}\\
-\Delta\twdnu=&-\ga\tVunu e^{\twunu}+\ga\tVdnu e^{\twdnu}
\end{split}
\right.
\end{equation}
in $\Onu$, where
$\tVunu(X)=V_{1,n}(\Yun+\eun X)$ and $\tVdnu(X)=V_{2,n}(\Yun+\eun X)$,
which is of the form \eqref{Liouvillen2} with $a=b=1$, $c=d=\ga$.
Moreover, \eqref{Vassumpt} is satisfied since $0\le \widetilde V_{i,n}^1\le C$ and
\[
\int_{\Onu}e^{\widetilde w_{i,n}^1}\,dX=\int_{\calO} e^{w_{i,n}}\,dX\le C.
\]
We also note that \eqref{wid} implies that 
\begin{equation}
\label{twuid}
\twunu(X)\le\twunu(0)=0,\qquad 
\twunu+\frac 1\ga \twdnu \to-\infty.
\end{equation}
We apply Lemma~\ref{lem:BMsystem}
on a ball $B_R\subset\mathbb R^2$ of fixed large radius $R>0$.
In view of \eqref{twuid} we rule out Alternative~1 and Alternative~3.
Suppose Alternative~2.\ holds. Then, $\twunu$ is uniformly bounded in $B_R$
and $\twdnu\to-\infty$ locally uniformly in $B_R$.
In particular, there exists $\twuu\in C_{\mathrm loc}^{1,\al}(\mathbb R^2)$ such that $\twunu\to\twuu$,
$e^{\twdnu}\to0$, locally uniformly in $\mathbb R^2$.
Consequently, we derive that $\twuu$ satisfies the Liouville equation
\[
-\Delta\twuu=V_1e^{\twuu}\qquad\mathrm{on\ }\mathbb R^2.
\]
Since
\[
\int_{B_R}e^{\twuu}=\lim_n\int_{B_R}e^{\twunu}\le\int_{\Onu}e^{\twunu}\le C
\]
we also have $\int_{\mathbb R^2}e^{\twuu}\le C$.
Now the classification theorem in \cite{cl} implies that
$\int_{\mathbb R^2}V_1e^{\twuu}=8\pi$ and consequently $m_1(p_0)=8\pi$.
This established the asserted lower bound $m_1(p_0)\ge8\pi$
in the case where $\twunu,\twdnu$ satisfy Alternative~2.
Suppose Alternative~4 holds.
Then $\twunu$ is locally uniformly bounded,$\mathcal S_1^0=\emptyset$,
$\mathcal S_2^0\neq\emptyset$ and $\twdnu\to-\infty$ locally uniformly.
Moreover,
\[
\widetilde V_{2,n}e^{\twdnu}\,dX\stackrel{\ast}{\rightharpoonup}m_2(x_0)\delta_0.
\]
We conclude that there exists $\twuu\in C_{\mathrm loc}^{1,\al}(\mathbb R^2)$
such that $\twunu\to\twuu$ locally uniformly, 
\[
-\Delta\twuu=V_1e^{\twuu}-m_2(x_0)\delta_{x_0}.
\]
In view of \cite{cl1} we conclude that
\[
\int_{\mathbb R^2}V_1e^{\twuu}>4\pi+m_2(x_0)>8\pi,
\]
and the asserted lower bound $m_1(x_0)\ge8\pi$ holds true
in this case as well.
\par
\textit{Proof of $m_2(p_0)\ge\frac{8\pi}{\ga}$.}
Similarly as above, we rescale \eqref{Liouvillen} around $\Ydn$ with respect to $\edn$ given by
\begin{align*}
\edn=&e^{-{w_{2,n}(\Ydn)} /2},
\end{align*}
that is $w_{2,n}(Y_{2,n})=\sup_{\calO}w_{2,n}=-2\log\edn$.
We define the expanding domain
\[
\Ond=\{X\in\mathbb R^2:\ \Ydn+\edn X\in\calO\}
\]
and we define $\twund,\twdnd:\Ond\to\mathbb R$ by setting
\begin{align*}
\twund(X)=&w_{1,n}(\Ydn+\edn X)-w_{2,n}(\Ydn)\\
\twdnd(X)=&w_{2,n}(\Ydn+\edn X)-w_{2,n}(\Ydn).
\end{align*}
Then, $\twund,\twdnd$ is a solution for the Liouville system
\begin{equation}
\label{Liouville2}
\left\{
\begin{split}
-\Delta\twund=&\tVund e^{\twund}-\tVdnd e^{\twdnd}\\
-\Delta\twdnd=&-\ga\tVund e^{\twund}+\ga\tVdnd e^{\twdnd}
\end{split}
\right.
\end{equation}
in $\Ond$, where $\tVund(X)=V_{1,n}(\Ydn+\edn X)$ and $\tVdnd(X)=V_{2,n}(\Ydn+\edn X)$.
Furthermore, as above,
\begin{equation*}
\label{twdid}
\twund+\frac 1\ga \twdnd \to-\infty.
\end{equation*}
We observe that
\[
0\leqslant \widetilde V_{i,n}^2(X)\leqslant C,
\qquad
\int_{\mathcal O_n^2}e^{\widetilde w_{i,n}^2}\,dX\leqslant  C
\]
for $i,=1,2$, for some $C>0$.
Therefore, Lemma~\ref{lem:BM}
may be applied locally to the Liouville systems \eqref{Liouville1} and \eqref{Liouville2}.
By analogous arguments as above, we conclude that there exists $\twdd$ such that $\twdnd\to\twdd$
locally uniformly in $\mathbb R^2$
with
\[
\ga\int_{\mathbb R^2}\widetilde V_2^2 e^{\twdd}\ge8\pi.
\]
That is,
\[
m_2(p_0)=\int_{\mathbb R^2}\widetilde V_2^2 e^{\twdd}\ge\frac{8\pi}{\ga},
\]
as desired.
%At this point, the remaining part of the proof of Proposition~\ref{prop:blowup}, Part~(i)  is completely analogous to
%\cite{OhtsukaSuzukiBanachCenter}.  We outline the main ideas for the reader's convenience. We use Lemma \ref{lem:BM}.
%n view of the identities \eqref{twuid}--\eqref{twdid}, we rule out Alternative~1.
%On the other hand, since $\twunu(0)=0=\twdnd(0)$, we rule out Alternative~3.
%Let $\widetilde w_k:\Rd\to\mathbb R$ be such that $\widetilde w_{k,n}^k\to\widetilde w_k$
%uniformly on compact subsets of $\Rd$, $k=1,2$.
%It follows that the limit function $\widetilde w_1$ of $\widetilde w_{1,n}^1$ solves 
%the standard Liouville equation
%\begin{equation*}
%-\Delta\tilde w_1 =V_1(0)e^{\tilde w_1} \qquad \mbox{ in }\Rd,\qquad \int_{\Rd} e^{\tilde w_1}<+\infty,
%\end{equation*}
%Hence, in view of \cite{cl}, we derive
%\[
%m_1\geqslant \int_{\Rd}V_1(0)e^{\tilde w_1}=8\pi.
%\]
\end{proof}
%%%%%%%%%%%%%%%%%%%%%%%%%%%%%%%%%%%%%%%%%%%%%%%%%%%%%%%%%%%%%%%%%%%%%%%%%%%%%%%%%%%%%%%%%%%%%%%%%%%%%%%%%%%%%%%%%%%%%%%%%%%%%%%%%%%%%
%%%%%%%%%%%%%%%%%%%%%%%%%%%%%%%%%%%%%%%%%%%%%%%%%%%%%%%%%%%%%%%%%%%%%%%%%%%%%%%%%%%%%%%%%%%%%%%%%%%%%%%%%%%%%%%%%%%%%%%%%%%%%%%%%%%%%
\section{Proof of Theorem~\ref{thm:mp}}
\label{sec:mp}
%%%%%%%%%%%%%%%%%%%%%%%%%%%%%%%%%%%%%%%%%%%%%%%%%%%%%%%%%%%%%%%%%%%%%%%%%%%%%%%%%%%%%%%%%%%%%%%%%%%%%%%%%%%%%%%%%%%%%%%%%%%%%%%%%%%%%
%%%%%%%%%%%%%%%%%%%%%%%%%%%%%%%%%%%%%%%%%%%%%%%%%%%%%%%%%%%%%%%%%%%%%%%%%%%%%%%%%%%%%%%%%%%%%%%%%%%%%%%%%%%%%%%%%%%%%%%%%%%%%%%%%%%%%
We recall from Section~\ref{sec:intro} that the variational functional for problem~\eqref{eq:mainpb}
is given by:
\[
J(u)=\frac12\int_\Sigma|\nabla_g v|^2-\lambda_1\log\int_\Sigma e^v-\frac{\lambda_2}{\gamma}\log\int_\Sigma e^{-\gamma v},
\quad v\in\mathcal E.
\]
We begin by checking that $J$ admits a mountain pass geometry for all $(\la_1,\la_2)\in\Lambda$.
\begin{lemma}[Existence of a local minimum]
\label{lem:localmin}
%Assume that $8\pi<{\mu _{1}(\Sigma)\left\vert \Sigma \right\vert}$ and let 
Fix $\lambda_1,\lambda_2 \geq 0.$ If 
\[
\lambda_1+\gamma \lambda_2 <{\mu _{1}(\Sigma )\left\vert \Sigma\right\vert},
\]
then $v\equiv0$ is a strict local minimum for $J$.
\end{lemma}
\begin{proof}
We set
\begin{equation*}
\label{g}
G_1(v):=\log\left( \int_{\Sigma }e^{ v}\,\dvg\right),\quad 
G_2(v):=\log\left( \int_{\Sigma }e^{ -\ga v}\,\dvg\right),
\end{equation*}
so that we may write
\[
J(v)=\frac{1}{2}\int_{\Sigma }\left\vert\nabla_g v\right\vert^{2}-\lambda_1G_1(v)-\frac{\lambda_2}{\ga}G_2(v).
\] 
For every $\phi\in \mathcal E$, we compute:
\begin{equation*}
G_1^{\prime}(v)\phi=\int_{\Sigma }\frac{e^{ v}\phi }{\int_{\Sigma }e^{v}};
\qquad G_2^{\prime}(v)\phi=\int_{\Sigma }\frac{-\ga e^{\ga v}\phi }{\int_{\Sigma}e^{-\ga v}}.
\end{equation*}
Moreover, for every  $\phi,\psi \in \mathcal E$,
\begin{equation*}
\label{g1sec}
\left\langle G_1^{\prime\prime}(v)\phi,\psi\right\rangle
=\frac{\left(\int_{\Sigma}e^{v}\phi\psi\right)
\left(\int_{\Sigma}e^{ v}\right) - \left( \int_{\Sigma}  e^{ v}\phi\right)
\left(\int_{\Sigma } e^{ v}\psi\right)}{\left(\int_{\Sigma } e^{ v}\right)^2 }.
\end{equation*}
\begin{equation*}
\label{g2}
\left\langle G_2^{\prime\prime}(v)\phi,\psi\right\rangle
=\ga^2 \frac{\left(\int_{\Sigma}e^{-\ga v}\phi\psi\right)
\left(\int_{\Sigma}e^{-\ga v}\right) - \left( \int_{\Sigma}  e^{-\ga v}\phi\right)
\left(\int_{\Sigma } e^{-\ga v}\psi\right)}{\left(\int_{\Sigma } e^{-\ga v}\right)^2 }.
\end{equation*}
In particular, for every $\phi \in \mathcal E$, we derive that
\begin{equation}
\begin{split}
G_1^{\prime}(0)\phi=G_2^{\prime}(0)\phi=0;\quad
\left\langle G_1^{\prime\prime}(0)\phi,\phi\right\rangle=\frac 1{|\Sigma|}\int_\Sigma \phi^2;\quad
\left\langle G_2^{\prime\prime}(0)\phi,\phi\right\rangle=\frac{\ga^2}{|\Sigma|}\int_\Sigma \phi^2.
\end{split}
\end{equation}
%Moreover,
%\[
%J^\prime (u)\phi=\int_\Sigma \left\langle \nabla v, \nabla \phi\right\rangle-\lambda_1 G_1^{\prime}(v)\phi- \frac {\lambda_\ga }\ga G_\ga^{\prime}(v)\phi.
%\]
Consequently, for every $ \phi\in \mathcal E$,
\begin{equation}\label{J1}
J^\prime(0)\phi=0, \qquad \qquad 
\left\langle J^{\prime\prime}(0)\phi,\phi\right\rangle=\int_\Sigma |\nabla \phi|^2
-\frac{{\lambda_1}+\ga \lambda_2}{|\Sigma|} \int_\Sigma \phi^2.
\end{equation}
Recalling the Poincar\'e inequality
\[
\int_\Sigma \phi^2 \leq \frac1{\mu_1(\Sigma)} \int_\Sigma |\nabla \phi|^2,\qquad \forall  \phi\in \mathcal E,
\]
where $\mu_1(\Sigma)$ is the first eigenvalue defined in \eqref{def:mu1},
we deduce by Taylor expansion at $0$ that
\begin{equation}\label{J2}
\begin{split}
J(u)-J(0)&= J^\prime (0)u +\frac  12  \left\langle J^{\prime\prime}(0)u,u\right\rangle+o(\|u\|^2)\\
&\geq \frac  12  \left(   1 -\frac{\lambda_1 + \ga \lambda_2 }{\mu_1(\Sigma) |\Sigma | } \right ) \| u\|+o (\| u\|^2),
\end{split}
\end{equation}
for all $u\in\mathcal E$.
Hence, $0$ is a local minimum for $J$, as asserted.
\end{proof}
In the next lemma we check that there exists a direction for $J$ along which $J$ is unbounded
from below.
\begin{lemma} 
\label{lem:Junbdd}
%The functional $J(v)$
%is bounded from below on $\mathcal E$ if and only if 
%\[
%\max \left\{  \lambda_1,\,\, \gamma  \lambda_2 \right\} \leq8\pi \\
%\]
If $\max\{ \lambda_1,\gamma \lambda_2 \}>8\pi$, then there exists $v_1\in\mathcal E$ such that 
$J(v_1)<0$ and $\|v_1\|\ge1$.
\end{lemma}
\begin{proof}  
The unboundedness from below of $J$, i.e., the existence of $v_1$, may be derived 
simply inserting the special
form of $\mathcal P$ given by \eqref{def:specialP} into formula~\eqref{def:barlambda}.
However, since the Struwe monotonicity argument also requires a control on $\|v_1\|$, we use a family of test functions
derived from the classical \lq\lq Liouville bubbles".
Namely, we fix $p_0\in\Sigma$ and $r_0>0$ a constant smaller than the injectivity radius of $\Sigma$
at $p_0$. Let $\mathcal B_{r_0}=\{p\in\Sigma:\ d_g(p,p_0)<r_0\}$ denote the geodesic ball
of radius $r_0$ centered at $p_0$.
For every $\eps>0$ let $\ue$ be the function defined by
\[
\ue(p)=
\begin{cases}
\log\frac{\eps^2}{(\eps^2+d_g(p,p_0)^2)^2}&\mbox{in\ }\mathcal B_{r_0}\\
\log\frac{\eps^2}{(\eps^2+r_0^2)^2}&\mbox{in\ }\Si\setminus\mathcal B_{r_0}.
\end{cases}
\]
Let $\ve\in\mathcal E$ be defined by
\[
\ve=\ue-\frac{1}{|\Si|}\int_\Si\ue\,\dvg.
\]
By elementary computations we have
\begin{equation*}
\begin{aligned}
&\int_\Si|\nabla_g\ue|^2\,\dvg=16\pi\log\frac{1}{\eps^2}+O(1);\\
&\frac{1}{|\Si|}\int_{\Si}\ue\,\dvg=\ln\eps^2+O(1);
\end{aligned}
\end{equation*}
and therefore
\begin{equation*}
\begin{aligned}
&\int_\Si|\nabla_g\ve|^2\,\dvg=16\pi\log\frac{1}{\eps^2}+O(1);\\
&\log\int_\Si e^{\ve}\,\dvg=\log\frac{1}{\eps^2}+O(1);\\
&\log\int_\Si e^{-a\ve}\,\dvg=O(1),\quad\forall a>0.
\end{aligned}
\end{equation*}
By assumption we have $\la_1>8\pi$ or $\la_2>8\pi/\ga$.
Suppose $\la_1>8\pi$.
Then, in view of the expansions above, we have:
\[
J(\ve)=(8\pi-\la_1)\log\frac{1}{\eps^2}+O(1),
\]
and therefore there exists $\eps_1>0$ such that $v_1:=v_{\eps_1}$
satisfies the desired properties.
Suppose $\la_2>8\pi/\ga$. Then,
\[
J(-\frac{\ve}{\ga})=\frac{1}{\ga}(\frac{8\pi}{\ga}-\la_2)\log\frac{1}{\eps^2}+O(1),
\]
and the statement holds true in this case as well.
\end{proof}
In particular,  Lemma~\ref{lem:localmin} and Lemma~\ref{lem:Junbdd} imply that 
the functional~$J$ admits a mountain-pass geometry whenever  $(\lambda_1,\lambda_2)\in\Lambda$.
\par
Finally, in order to conclude the proof of Theorem~\ref{thm:mp} we need the
following consequences of Theorem~\ref{thm:blowup}.
\begin{lemma}
\label{lem:blowup}
The following properties hold.
%
%  \item [(i)]
%{\em Improved minimal mass:}
%The lower bound~\eqref{8pi} is improved as follows:
%\begin{equation}
%\label{8pi2}
%\mm(p)\ge8\pi\ \mbox{for all $p\in\Sp$},\qquad\qquad \mmm(p)\ge\frac{8\pi}{\ga}\ \mbox{for all $p\in\Sm$}.
%\end{equation}
\begin{enumerate}\label{thm:improvbu}
\item[(i)]
For every $p\in\Sp\cap\Sm$ we have
\begin{equation}
\label{improvedsum}
\mm+\ga\mmm\ge 16\pi(1+{\ga}).
\end{equation}
\item[(ii)]
If
$\Sp\setminus\Sm\neq\emptyset$, then $r_1\equiv 0$; in particular, in this case $\lambda_{1,0}= 8\pi\,\sharp\Sp$.
Similarly, if $\Sm\setminus\Sp\neq\emptyset$, then $r_2\equiv0$ and $\la_{2,0}= \frac{8\pi}\ga\,\sharp\Sm$.
\end{enumerate}
\end{lemma}
\begin{proof}
Proof of Part~(i).
In view of Theorem~\ref{thm:blowup} and identity~\eqref{identquad} it suffices to show that
\begin{equation}
\al_\ga:=\min\left\{x+\ga y:\ x\ge8\pi,\ y\ge\frac{8\pi}{\ga},\ (x-y)^2=8\pi(x+\frac{y}{\ga})\right\}
=16\pi(1+\ga).
\end{equation}
Let $\mathbf p:=\{(x,y)\in\mathbb R^2:\ (x-y)^2=8\pi(x+\frac{y}{\ga})\}$
be the parabola defined by the quadratic relation \eqref{identquad}.
Let $\overline x,\overline y>0$ be such that $(\overline x,\frac{8\pi}{\ga})\in\mathbf p$
and $(8\pi,\overline y)\in\mathbf p$.
In view of the geometric properties of $\mathbf p$ we find that
\[
\al_\ga=\min\{\overline x+8\pi,8\pi+\ga\overline y\}.
\] 
Solving the equation $(x-y)^2=8\pi(x+\frac{y}{\ga})$ with respect to $x$, we have
\begin{align*}
x(y)=x_\pm(y)=y+4\pi\pm\sqrt{8\pi(1+\frac{1}{\ga})y+16\pi^2}.
\end{align*} 
We deduce that
\begin{align*}
\overline x=x_+(\frac{8\pi}{\ga})=8\pi(1+\frac{2}{\ga}).
\end{align*}
Similarly, solving $(x-y)^2=8\pi(x+\frac{y}{\ga})$ with respect to $y$,
we have
\[
y=y_\pm(x)=x+\frac{4\pi}{\ga}\pm\sqrt{8\pi(1+\frac{1}{\ga}x+\frac{16\pi^2}{\ga^2})}.
\]
We conclude that
\[
\overline y=y_+(8\pi)=8\pi(2+\frac{1}{\ga}).
\]
It follows that
\[
\al_\ga=\min\{16\pi(1+\frac{1}{\ga}),\,16\pi(1+\ga)\}=16\pi(1+\ga).
\]
Proof of Part~(ii).
Let
\begin{align*}
&\uun=G\ast\llan\frac{e^{\vn}}{\int_\Si e^{\vn}},
&&\udn=G\ast\lllan\frac{e^{-\ga\vn}}{\int_\Si e^{-\ga\vn}},
\end{align*}
where $G=G(x,p)$ denotes the Green's function defined by
\begin{equation*}
\begin{cases}
-\Delta_{x}G(\cdot,p)=\delta_p-\frac{1}{|\Si|}&\mbox{on\ }\Si\\
\int_\Si G(\cdot,p)\,\dvg=0
\end{cases}
\end{equation*}
and where $\ast$ denotes convolution.
Then, $\vn=\uun-\udn$.
For a measurable function $f$ defined on a subset of $\Si$ and
for any $T>0$ we set $f^T=\min\{f,T\}$.
Suppose $p\in\Sp\setminus\Sm$. 
We claim that
\begin{equation}
\label{eq:intinfty}
\int_\Si e^{\vn}\to+\infty.
\end{equation}
Indeed, we have, using \eqref{eq:concentration} and Theorem~\ref{thm:blowup}:
\begin{align*}
\uun(x)\ge G^T\ast\llan\frac{e^{\vn}}{\int_\Si e^{\vn}}\to G^T\ast\sum_{p\in\Sp}(\mm(p)+r_1)
\ge8\pi G^T(x,p)-C.
\end{align*}
In local coordinates centered at $p$ we obtain
\[
\uun(X)\ge\ln[\frac{1}{|X|^4}]^T-C.
\]
Since $p\not\in\Sm$, we have that $\|\udn\|_{L^\infty(B_\rho(p))}$
is uniformly bounded with respect to $n$, for some sufficiently small $\rho>0$.
Now, Fatou's lemma implies that
\[
\liminf_{n\to\infty}\int_\Si e^{\vn}
=\liminf_{n\to\infty}\int_\Si e^{\uun-\udn}
\ge c_0\int_{B_\rho(p)}[\frac{1}{|X|^4}]^T
\]
for some $c_0>0$ independent of $n$.
Letting $T\to+\infty$, we deduce that \eqref{eq:intinfty} holds true,
and consequently $r_1\equiv0$.
\par
Now, we assume that $p\in\Sm\setminus\Sp$. Similarly as above, we claim that
\begin{equation}
\label{eq:intgainfty}
\int_\Si e^{-\ga\vn}\to+\infty.
\end{equation}
Indeed, we note that
\begin{align*}
\udn(x)\ge G^T\ast\lllan\frac{e^{-\ga\vn}}{\int_\Si e^{-\ga\vn}}
\to G^T\ast\sum_{p\in\Sm}(\mmm(p)\delta_p+r_2)
\ge\frac{8\pi}{\ga}\sum_{p\in\Sm}G^T(x,p)-C,
\end{align*}
for some $C>0$ independent of $n$.
In particular, in local coordinates $X$ near $p$,
\[
\udn(X)\ge\frac{4}{\ga}\ln[\frac{1}{|X|}]^T-C.
\]
On the other hand, since $p\not\in\Sp$,
we have that $\|\uun\|_{L^\infty(B_\rho(p))}$ is bounded independently of $n$
in a ball centered at $p$ provided $\rho>0$ is sufficiently small.
We conclude that
\[
\int_\Si e^{-\ga\vn}\ge c_0\int_{B_\rho(p)}[\frac{1}{|X|^4}]^T.
\]
In view of Fatou's lemma and letting $T\to+\infty$, we derive \eqref{eq:intgainfty}.
Therefore, $r_2\equiv0$, as asserted.
\end{proof}
Finally, we provide the proof of Theorem \ref{thm:mp}.
\begin{proof}[Proof of Theorem \ref{thm:mp}]
In view of Struwe's Monotonicity Trick, for almost every $(\lambda_1,\lambda_2) \in \Lambda$ there exists a nontrivial solution to \eqref{eq:mainpb}.
As this argument is by now standard, we refer to \cite{Struwe, StTa} for a proof. A detailed proof in this
specific context may also be found in \cite{RZ}.
Let $(\lambda_1^0,\lambda_2^0)\in\Lambda$ be fixed. Let  $(\lambda_{1,n},\la_{2,n},v_n)$ be a solution sequence 
for \eqref{eq:mainpb} with $\la_{1,n}\to\la_{1}^0$ and $\la_{2,n}\to\la_{2}^0$, 
$(\lambda_{1,n},\la_{2,n})\in\Lambda$. We claim that $v_n$ converges to some $v\in \mathcal E$, 
solution to \eqref{eq:mainpb}. To this aim we exploit the blow up analysis as derived 
in Theorem \ref{thm:blowup}. By contradiction we assume that $p_0\in\mathcal S$.  
%ote that in view of \eqref{def:Lambda}, $\mathcal S$ consists at most in one blow up point $p_0\in \Sigma$. 
%hence, we may assume by contradiction that  $S=\{p_0\}$. 
If $p_0\in\mathcal S_1\setminus\Sm$, then in view of Lemma~\ref{lem:blowup}--(ii)
we have $\la_1^0\in8\pi\mathbb N$, a contradiction to \eqref{def:Lambda}.
If $p_0\in\Sm\setminus\Sp$, then in view of Lemma~\ref{lem:blowup}--(ii)
we have $\la_2^0\in\frac{8\pi}{\ga}\mathbb N$, a contradiction to \eqref{def:Lambda}.
Suppose that $p_0\in \mathcal S_1\cap\mathcal S_2$. In this case, 
using  \eqref{improvedsum},
\[
\lambda_1+\gamma \lambda_2  \geq m_1(p_0)+\gamma m_2 (p_0)\geq 16 \pi (1+\gamma)
\]
and we again obtain a contradiction to \eqref{def:Lambda}. 
%Hence $\mathcal S_1\cap \mathcal S_2=\emptyset.$  
%Suppose $p_0\in \mathcal S_1= \mathcal S_1\setminus \mathcal S_2.$ 
%In this case $\lambda_1=8\pi$ and then, by using  \eqref{def:Lambda} it is  $\ga\lambda_2>8\pi$, a contradiction. 
%Hence, $\mathcal S_1=\emptyset.$  Finally, suppose that $p_0\in \mathcal S_2=\mathcal S_2\setminus \mathcal S_1$. 
%In this case  $\lambda_2=\frac {8\pi}\ga$ and so  $\lambda_1>8\pi$,  a contradiction. We deduce that $\mathcal S=\emptyset$. 
Hence, the compactness for the solution sequence  $(\lambda_{1,n},\la_{2,n},v_n)$ holds and Theorem~\ref{thm:mp} 
is completely established.
\end{proof}
{\bf Acknowledgments}.
This work is supported by Progetto GNAMPA-INDAM 2015: \emph{Alcuni aspetti di equazioni ellittiche non lineari}.
The first author is also supported by PRIN {\em Aspetti variazionali e perturbativi nei problemi differenziali nonlineari}.
The second author is also supported by Progetto STAR 2015: \emph{Variational Analysis and Equilibrium Models in Physical
and Socio-Economic Phenomena.}
%%%%%%%%%%%%%%%%%%%%%%%%%%%%%%%%%%%%%%%%%%%%%%%%%%%%%%%%%%%%%%%%%%%%%%%%%%%%%%%%%%%%%%%%%%%%%%%%%%%%%%
%%%%%%%%%%%%%%%%%%%%%%%%%%%%%%%%%%%%%%%%%%%%%%%%%%%%%%%%%%%%%%%%%%%%%%%%%%%%%%%%%%%%%%%%%%%%%%%%%%%%%%
%%%%%%%%%%%%%%%%%%%%%%%%%%%%%%%%%%%%%%%%%%%%%%%%%%%%%%%%%%%%%%%%%%%%%%%%%%%%%%%%%%%%%%%%%%%%%%%%%%%%%%
%%%%%%%%%%%%%%%%%%%%%%%%%%%%%%%%%%%%%%%%%%%%%%%%%%%%%%%%%%%%%%%%%%%%%%%%%%%%%%%%%%%%%%%%%%%%%%%%%%%%%%

\end{document}